\documentclass[report]{owrart}

\usepackage{array}
\usepackage{graphicx}
\graphicspath{{./Figures/}}

\def\bfp{\boldsymbol{p}}
\def\bfv{\boldsymbol{v}}

\def\bfxi{{\boldsymbol{\xi}}}

\def\BBB{\mathcal{B}}
\def\SSS{\mathcal{S}}

\subjclass{}

\newcommand{\SimS}[1]{\raisebox{-0.75em}{\includegraphics[scale=0.24]{SimplexSplineSmall#1.pdf}}}

\usepackage{tikz}

\def\PS{
\begin{tikzpicture}[scale = 1.7]
\draw[line width = 0.5] (-0.1, 0.0) -- (0.1,0.0);
\draw[line width = 0.5] (-0.1, 0.0) -- (0.0,0.1732);
\draw[line width = 0.5] ( 0.1, 0.0) -- (0.0,0.1732);
\end{tikzpicture}
}

\def\PSB{
\begin{tikzpicture}[scale = 1.7]
\draw[line width = 0.5] (-0.1, 0.0) -- (0.1,0.0);
\draw[line width = 0.5] (-0.1, 0.0) -- (0.0,0.1732);
\draw[line width = 0.5] ( 0.1, 0.0) -- (0.0,0.1732);

\draw[very thin]  (-0.05, 0.5*0.1732) -- (0.05, 0.5*0.1732);
\draw[very thin]  (-0.05, 0.5*0.1732) -- (0.0, 0.0);
\draw[very thin]  (0.0, 0.0) -- (0.05, 0.5*0.1732);

\draw[very thin]  (0.0, 0.0) -- (0.0,0.1732);
\draw[very thin] (-0.1, 0.0) -- (0.05, 0.5*0.1732);
\draw[very thin] (-0.05, 0.5*0.1732) -- (0.1,0.0);
\end{tikzpicture}
}

\newtheorem{theorem}{Theorem}

\makeatletter
\newenvironment{my_abstract}{%
    \small
    \quotation{\bfseries Abstract.}}
    {\endquotation}
    
\makeatother

\begin{document}
\begin{talk}{Tom Lyche and Georg Muntingh}
{Simplex Spline Bases on the Powell-Sabin 12-Split: Part II}
{Muntingh, Georg}

\begin{my_abstract}
For the space $\mathcal{S}$ of $C^3$ quintics on the Powell-Sabin 12-split of a triangle, we determine the simplex splines in $\mathcal{S}$ and the six symmetric simplex spline bases that reduce to a B-spline basis on each edge, have a positive partition of unity, a (barycentric) Marsden identity, and domain points with an intuitive control net. We provide a quasi-interpolant with approximation order 6 and a Lagrange interpolant at the domain points. The latter can be used to show that each basis is stable in the $L_\infty$ norm, which yields an $h^2$ bound for the distance between the B\'ezier ordinates and the values of the spline at the corresponding domain points. Finally, for one of these bases we provide $C^0$, $C^1$, and $C^2$ conditions on the control points of two splines on adjacent macrotriangles, and a conversion to the Hermite nodal basis.
\end{my_abstract}

\bigskip

\noindent Analogous to the $C^1$ quadratic simplex spline basis from \cite{GM_CohenLycheRiesenfeld13}, we derive $C^3$ quintic simplex spline bases on the Powell-Sabin 12-split $\PSB$ of a triangle $\PS$ \cite{GM_LycheMuntingh15}. The resulting computations are implemented in a {\tt Sage} worksheet, which can be downloaded and tried out online in {\tt SageMathCloud} \cite{GM_WebsiteGeorg}. We follow the notation in Part I.

A case-by-case analysis of the possible knot multiplicities yields:
\begin{theorem}\label{thm:AllowedSimplexSplines}
With one representative for each $S_3$ equivalence class, these are the $C^3$ quintic simplex splines on $\PSB$ that reduce to a B-spline on the boundary of $\PS$:
\begin{center}
\begin{tabular}{cccccccccc}
   \SimS{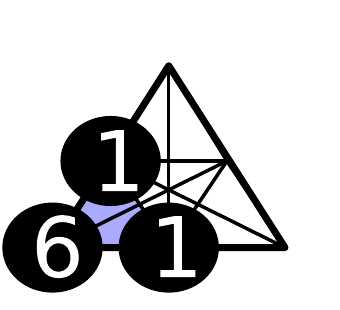} & \SimS{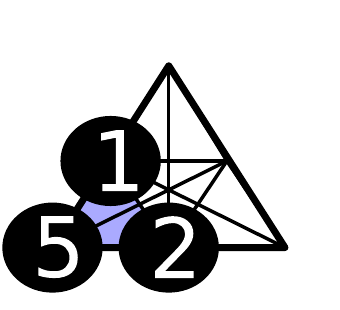} & \SimS{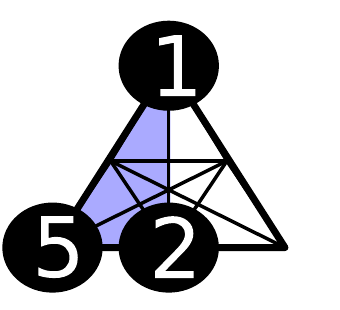} & \SimS{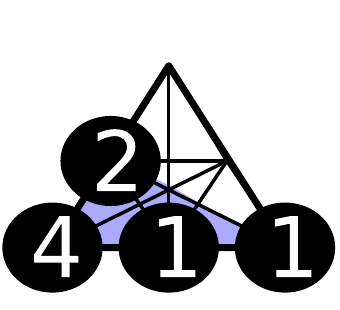} & \SimS{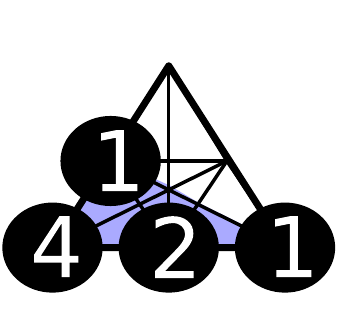}
 & \SimS{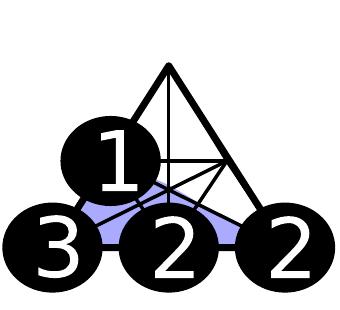} & \SimS{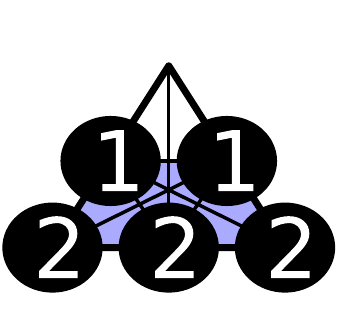} & \SimS{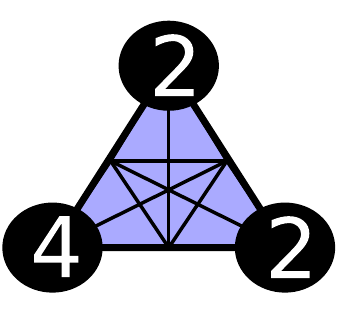} & \SimS{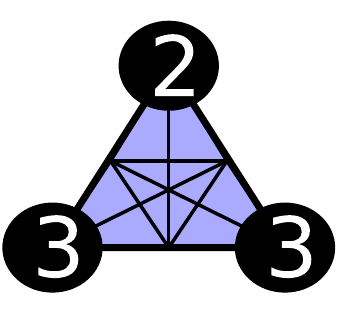} & \SimS{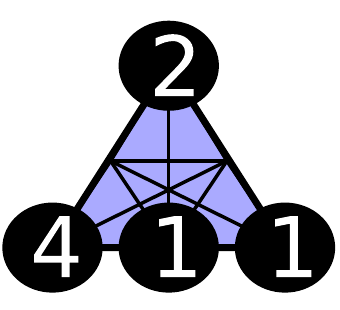}\\
   \SimS{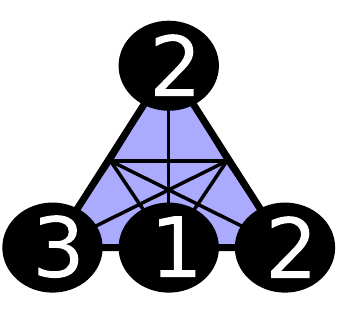} & \SimS{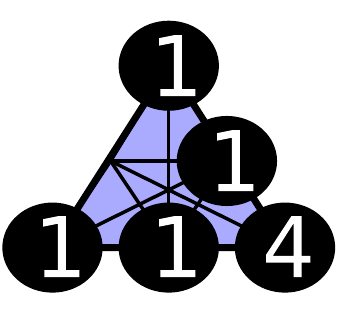} & \SimS{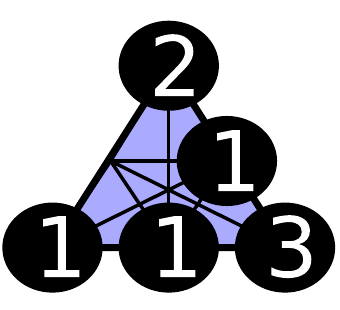} & \SimS{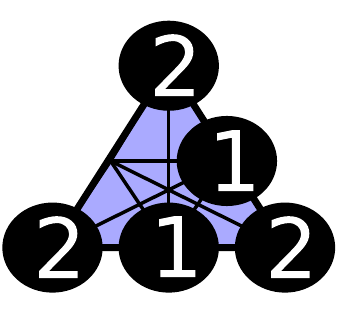} & \SimS{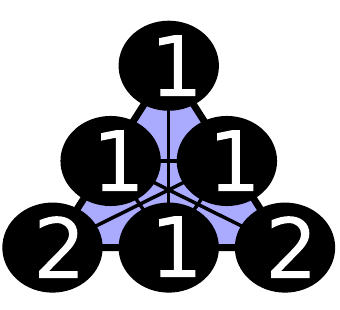}
 & \SimS{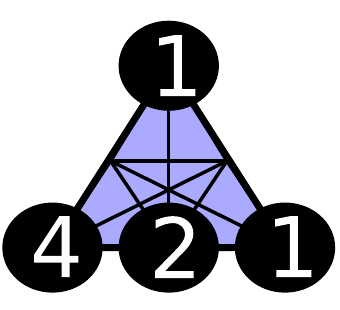} & \SimS{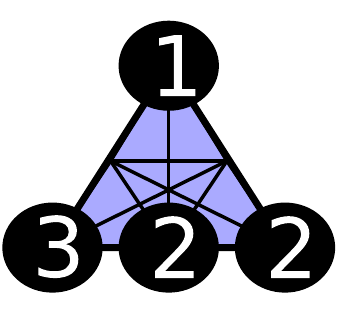} & \SimS{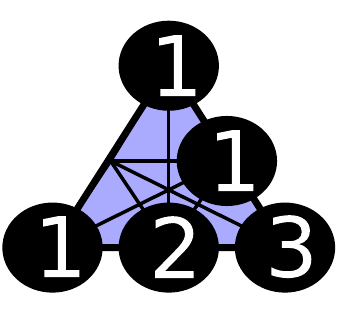} & \SimS{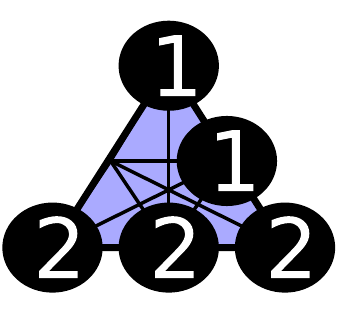} & \SimS{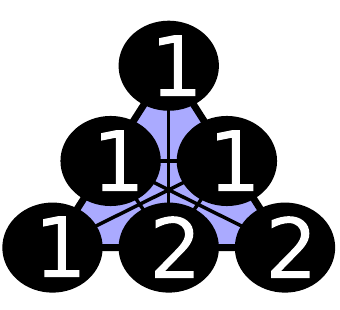}
\end{tabular}
\end{center}
\end{theorem}

We first create a large list of potential bases for the space $\SSS^3_5(\PSB)$ of $C^3$ quintics on the 12-split. Using the macro-element from \cite{GM_LycheMuntingh14}, we then narrow this down to a short list with good properties:

\begin{theorem}
There are precisely six sets $\BBB = \BBB_a, \BBB_b, \BBB_c,\BBB_d,\BBB_e, \BBB_f$ satisfying:
\begin{enumerate}
\item $\BBB$ is a basis of $\SSS^3_5(\PSB)$ consisting of simplex splines.
\item $\BBB$ is $S_3$-invariant.
\item $\BBB$ reduces to a B-spline basis on the boundary.
\item $\BBB$ has a positive partition of unity and a Marsden identity, for which the dual polynomials have only real linear factors.
\item $\BBB$ has all its domain points inside the macro triangle $\PS$, with precisely 8 domain points on each edge of $\PS$.
\end{enumerate}
\end{theorem}
For instance, the basis $\BBB_c = \{S_j\}_{j=1}^{39}$ is
\[ \left[ \frac14 \SimS{600101}, \frac14 \SimS{500201}, \frac12 \SimS{410201}, \frac12 \SimS{320201}, \frac34 \SimS{220211}, \SimS{141110}, \frac12 \SimS{131210}, \frac34 \SimS{121211} \right]_{S_3} \]
and satisfies the barycentric Marsden identity
\begin{align*}
& \qquad\qquad\qquad\qquad\qquad\left(c_1\SimS{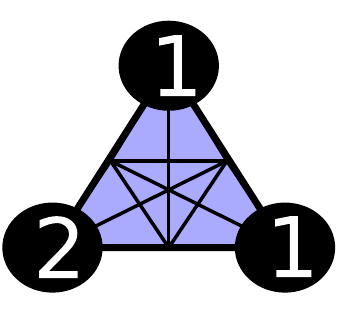} + c_2\SimS{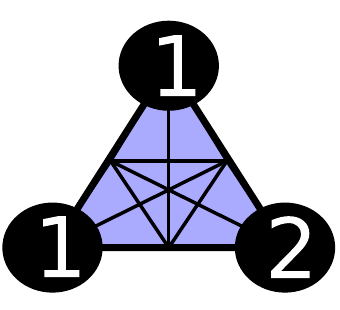} + c_3\SimS{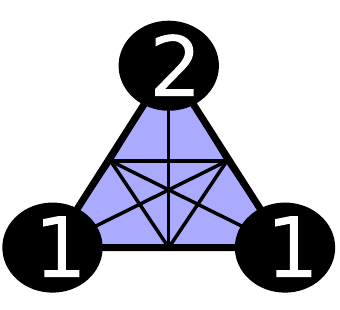} \right)^5 =\\
& \sum \left[\frac14 c_1^5 \SimS{600101}\right]_{S_3} \cup \left[\frac14 c_1^4 c_4 \SimS{500201}\right]_{S_3} \cup \left[\frac12 c_1^2 c_2 c_4^2 \SimS{320201}\right]_{S_3} \cup \left[\frac34 c_1 c_2 c_4 c_5 c_{10} \SimS{121211}\right]_{S_3}\\
& \cup \left[\frac12 c_1^3 c_4^2\SimS{410201}\right]_{S_3}\cup \left[\frac12 c_1 c_2^2 c_4 c_5 \SimS{131210}\right]_{S_3} \cup \left[\frac34 c_1c_2c_4^2c_{10} \SimS{220211}\right]_{S_3} \cup \left[c_2^3 c_4 c_5 \SimS{141110}\right]_{S_3}
\end{align*}

Factoring the dual polynomials and replacing `$c_i$' by `$\bfv_i$', one obtains 39 sets $\{\bfp_{j,r}^*\}_{r=1}^5$, $j=1,\ldots,39$, of dual points. Taking the average of each set one arrives at the domain points $\{\bfxi_j\}_{j=1}^{39}$. To preserve the symmetry of $\PSB$, the domain points are forced to form a hybrid mesh with triangles, quadrilaterals, and a hexagon in the center. This mesh is shown below on two adjacent macro triangles, together with an ordering of the domain points.

\begin{center}\includegraphics[scale=0.67, clip = true, trim = 7 8 45 35]{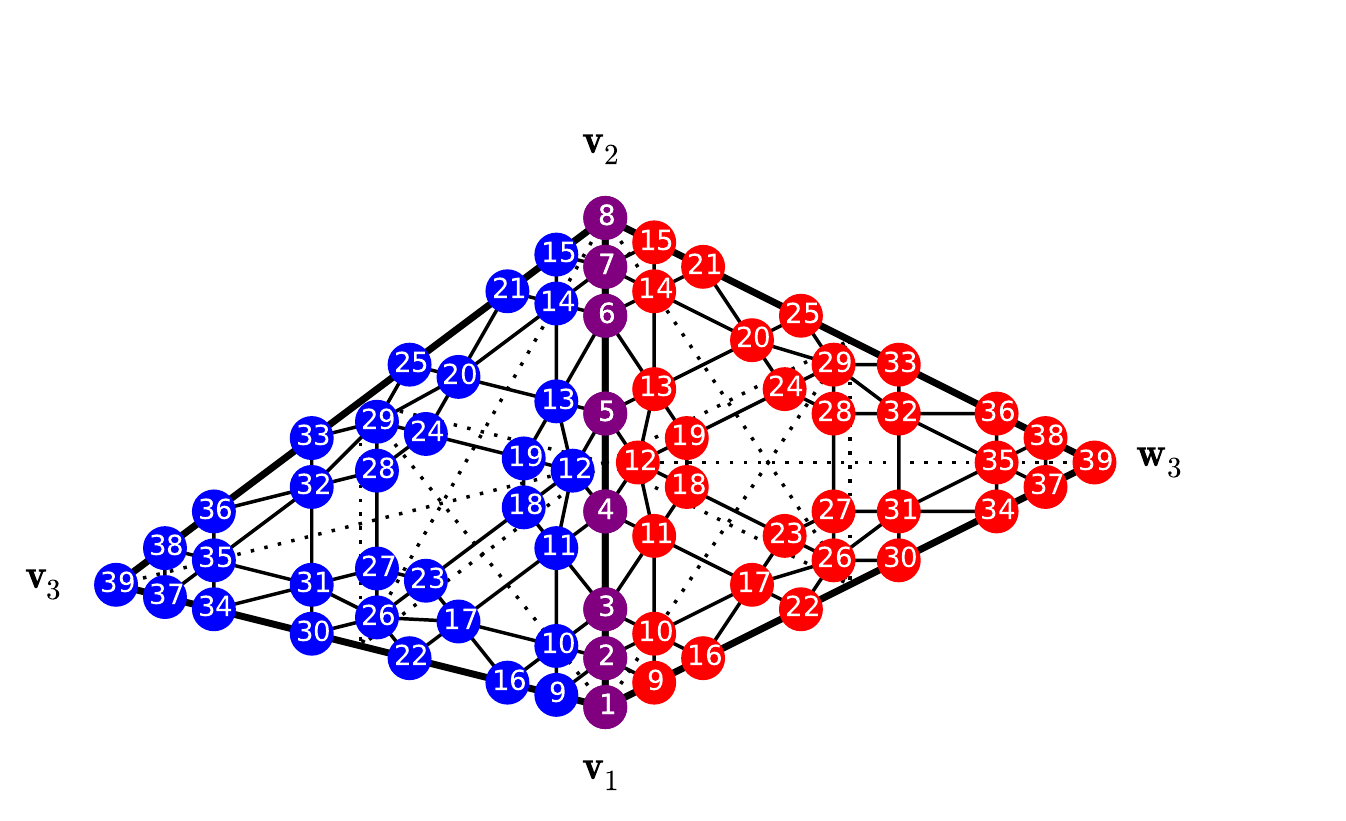}\end{center}

The collocation matrix $\{S_j(\bfxi_i\}_{i,j=1}^{39}$ is nonsingular, showing that $\{\bfxi_i\}_{i=1}^{39}$ is unisolvent for $\SSS^3_5(\PSB)$, i.e., there is a unique Lagrange interpolant at the domain points. Moreover, it was previously shown that there is a unique Hermite interpolant for the space $\SSS^3_5(\PSB)$ based on values and derivatives at the corners, midpoints, and quarterpoints \cite{GM_LycheMuntingh14}. Finally, the Marsden identity yields that
\[
Q(f) := \sum_{j=1}^{39} S_j \sum_{k=1}^5 \frac{1}{5!} k^5 (-1)^{k-1} \sum_{1\leq r_1 < \cdots < r_k \leq 5} f\left(\frac{\bfp_{j,r_1}^* + \cdots  + \bfp_{j,r_k}^*}{k}\right)
\]
is a quasi-interpolant that reproduces all polynomials up to degree 5 and has approximation order $6$. 
Moreover, using the Lagrange interpolant we show that the six bases are stable in the $L_\infty$ norm with a condition number bounded independently of the geometry.
As a consequence we obtain an $h^2$ bound of the distance between the B\'ezier ordinates and the values of the spline at the corresponding domain points. 

As in the above figure, let $\PS := [\bfv_1, \bfv_2, \bfv_3]$ and $\tilde{\PS} := [\bfv_1, \bfv_2, \tilde{\bfv}_3]$ be triangles sharing the edge $[\bfv_1, \bfv_2]$. Imposing a smooth join along $[\bfv_1, \bfv_2]$ of 
\[ f(\bfv) := \sum_{i=1}^{39} c_i S_i (\bfv), \ \bfv\in \PS,\qquad
\tilde{f}(\bfv) := \sum_{i=1}^{39} \tilde{c}_i \tilde{S}_i (\bfv), \ \bfv\in \tilde{\PS} \]
translates into linear relations among the B\'ezier ordinates $c_i$ and $\tilde{c}_i$.

\begin{theorem}
Let $(\beta_1,\beta_2,\beta_3)$ be the barycentric coordinates of $\tilde{\bfv}_3$ with respect to the triangle $\PS$. Then $f$ and $\tilde{f}$ meet with 

\noindent$C^0$ smoothness if and only if $\tilde{c}_i = c_i$, for $i = 1,\ldots,8$;\\
\noindent$C^1$ smoothness if and only if in addition\\
~\qquad\begin{tabular}{lllll}
&$\tilde{c}_9    = \beta_1c_1 + \beta_2c_2 + \beta_3c_9   $, &
&$\tilde{c}_{11} = \beta_1(2c_3 - c_2) + \beta_2c_4 + \beta_3c_{11}$, \\
&$\tilde{c}_{10} = \beta_1c_2 + \beta_2c_3 + \beta_3c_{10}$, &
&$\tilde{c}_{12} = \beta_1\frac{2c_4 + c_5}{3} + \beta_2\frac{c_4 + 2c_5}{3} + \beta_3 c_{12}$,
\end{tabular}\\
and analogous conditions for $\tilde{c}_{13}, \tilde{c}_{14}$, and $\tilde{c}_{15}$;

\noindent$C^2$ smoothness if and only if in addition \\
\begin{tabular}{lll}
& $\tilde{c}_{16} = \beta_1^2 c_1 + 2\beta_1\beta_2 c_2 + \beta_2^2c_3 + 2\beta_1\beta_3 c_9 + 2\beta_2\beta_3 c_{10} + \beta_3^2 c_{16}$,\\

& $\tilde{c}_{17} = \beta_1^2 c_2 + \beta_2^2 c_4 + \beta_3^2c_{17}+ 2\beta_1\beta_2 \frac{3c_3 - c_2}{2} + 2\beta_1\beta_3 \frac{3c_{10} - c_2}{2} + 2\beta_2\beta_3 \frac{c_{10} + 2c_{11} - c_3}{2}$,\\

& $\tilde{c}_{18} = \beta_1^2 \frac{2c_3 + 2c_4 - c_2 }{3} + \beta_2^2 \frac{c_4 + 2c_5}{3} + \beta_3^2c_{18}
+ 2\beta_1\beta_2 \frac{c_2 - 2c_3 + 6 c_4 + c_5}{6} $\\
& \qquad $ + 2\beta_1\beta_3 \frac{c_2 - 2c_3 + 2c_4 - c_5 + 3c_{11} + 3 c_{12}}{6} + 2\beta_2\beta_3 \frac{9c_{12} - 2c_5 - c_{11}}{6}$, \\
\end{tabular}\\
and analogous conditions for $\tilde{c}_{19}, \tilde{c}_{20}$, and $\tilde{c}_{21}$.
\end{theorem}
Whenever the domain points follow the shape of the macro triangles, we recover the classical B\'ezier conditions. All conditions are valid for the domain points as well, so that they also hold for the control points. Although conditions for $C^3$ smoothness can also be derived, one of these involves only $(\beta_1,\beta_2,\beta_3)$ and the B\'ezier ordinates on one triangle, showing that this element cannot be used to obtain $C^3$ smoothness on a general triangulation.

One can easily convert between $\BBB_c$ and the Hermite nodal basis from \cite{GM_LycheMuntingh14}. For instance, the nodal function corresponding to the point evaluation at $\bfv_1$ is
\begin{align*}
\varepsilon_{\bfv_1}^*    &  = \frac14 \SimS{600101} +  \frac14 \left(\SimS{500201} + \SimS{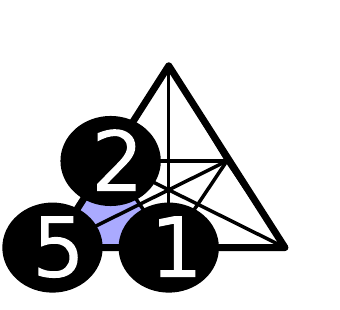}\right) + \frac12 \left(\SimS{410201} + \SimS{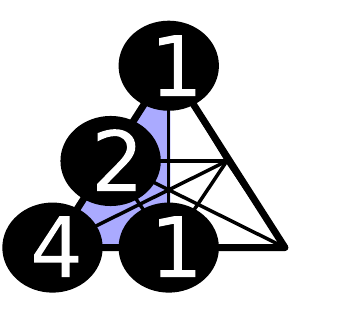}\right)  +  \SimS{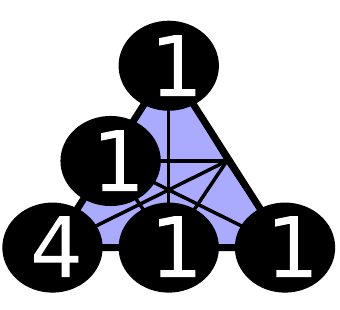}\\
&  +  \frac12 \left(\SimS{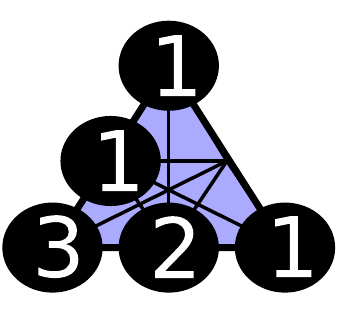} + \SimS{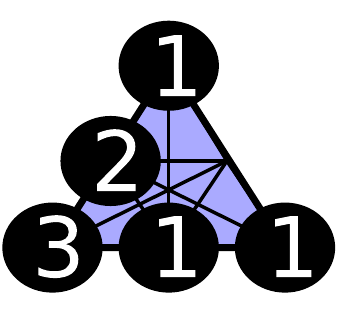}\right) + \frac12 \left(\SimS{320201} + \SimS{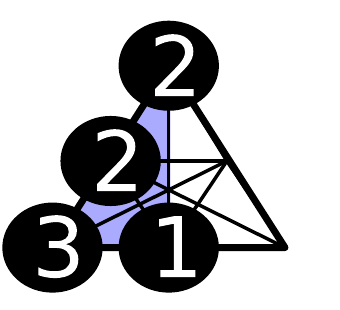}\right)+ \frac{9}{16} \left(\SimS{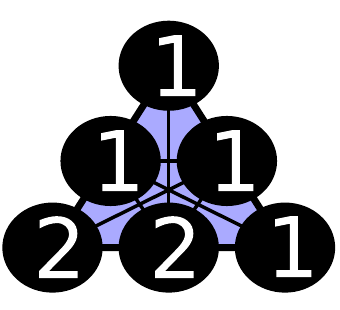} + \SimS{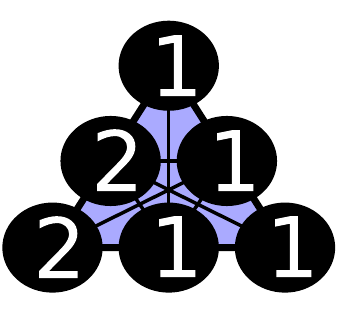}\right) \\
&  + \frac38 \left(\SimS{220211} + \SimS{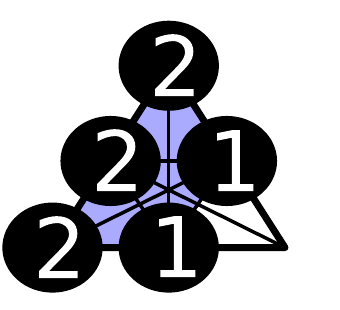}\right)
 + \frac{3}{16} \left(\SimS{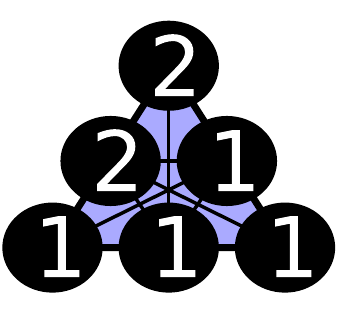} + \SimS{121211}\right)
\end{align*}
which, on a regular hexagon split at its barycenter, has the graph and wireframe
\begin{center}
\includegraphics[scale = 0.125, clip = true, trim = 120 8 125 155]{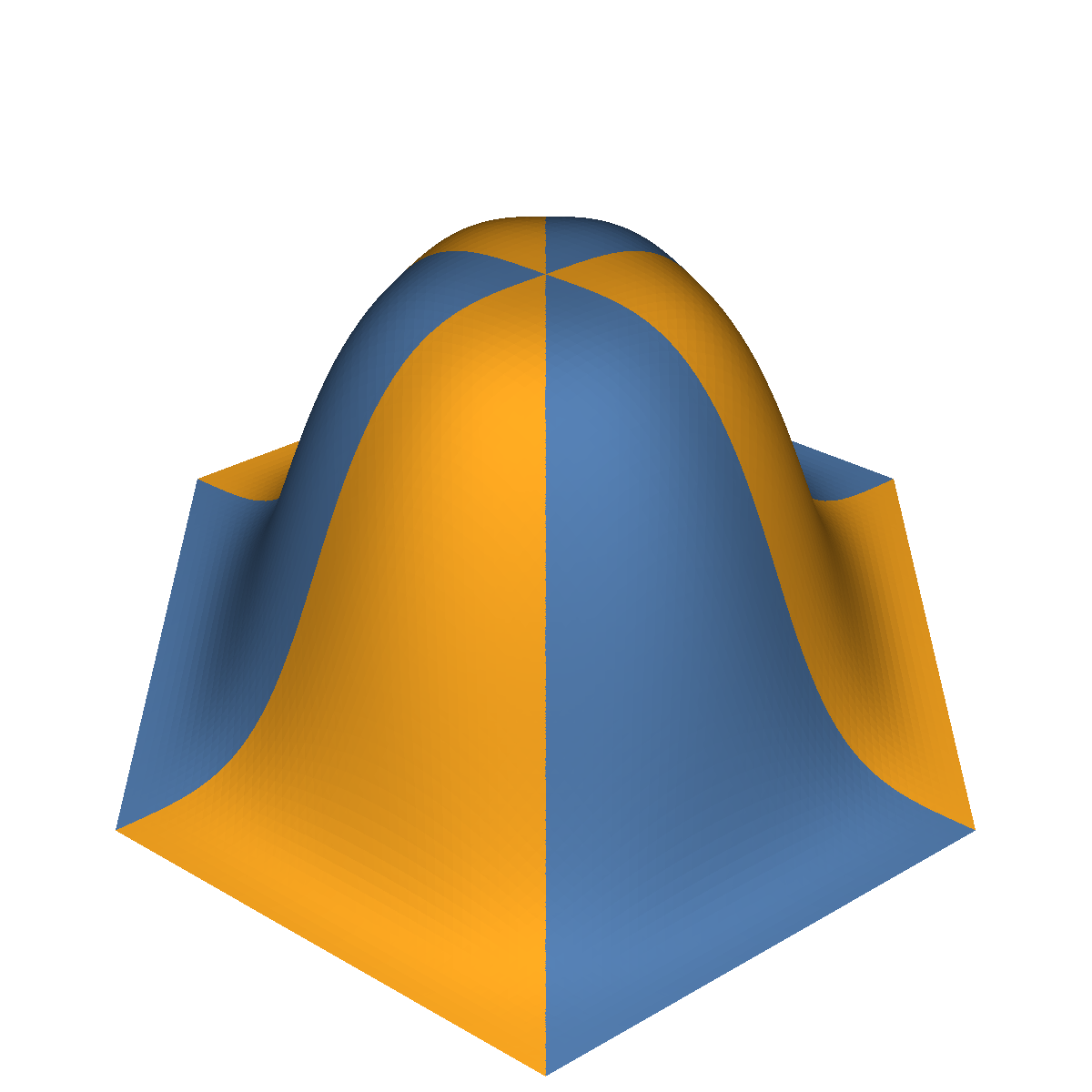}\qquad
\includegraphics[scale = 0.125, clip = true, trim = 120 8 125 155]{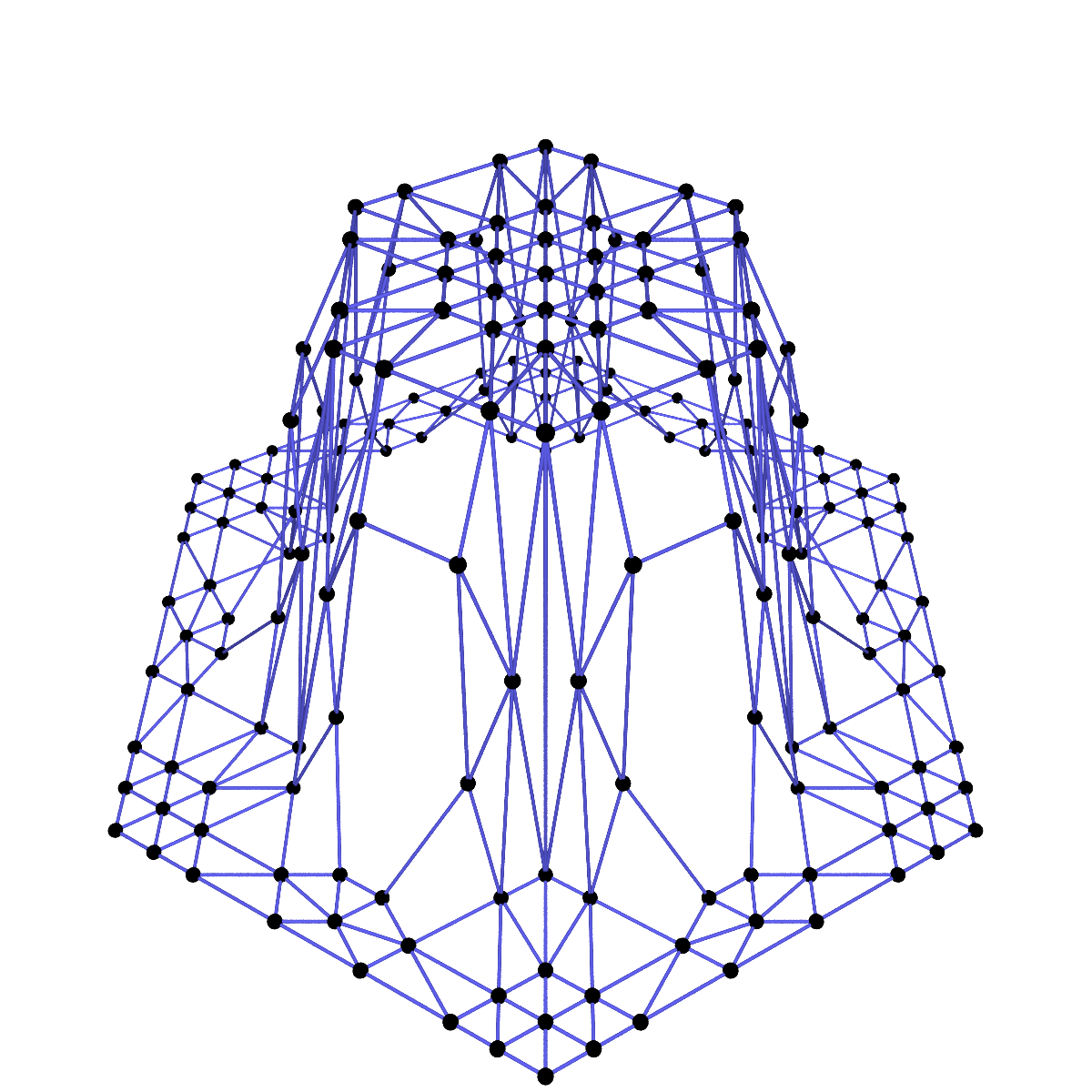}\\
\end{center}

\end{talk}


\begin{thebibliography}{99}

\bibitem{GM_CohenLycheRiesenfeld13}
Elaine Cohen, Tom Lyche, Richard Riesenfeld, \textit{A B-spline-like basis for the Powell-Sabin 12-split based on simplex splines}, Math. Comp. \textbf{82} (2013), no. 283, 1667--1707.

\bibitem{GM_LycheMuntingh14}
Tom Lyche and Georg Muntingh, \textit{A Hermite interpolatory subdivision scheme for $C^2$-quintics on the Powell-Sabin 12-split}, Comput. Aided Geom. Design \textbf{31} (2014), no. 7--8, 464--474.

\bibitem{GM_LycheMuntingh15}
Tom Lyche and Georg Muntingh, \textit{Stable simplex spline bases for $C^3$ quintics on the Powell-Sabin 12-split}, Available at \texttt{http://arxiv.org/abs/1504.02628}.

\bibitem{GM_WebsiteGeorg}
Georg Muntingh, \textit{Personal Website},\\ \texttt{https://sites.google.com/site/georgmuntingh/academics/software}.

\end{thebibliography}
\end{document}